\documentclass[a4paper, 11pt]{amsart}
\usepackage{amsmath, amssymb, eucal, amscd, amstext, enumerate, mathrsfs, yfonts, amsfonts}
\usepackage[plainpages=false,colorlinks,hyperindex,pdfpagemode=None,bookmarksopen,linkcolor=red,citecolor=blue,urlcolor=blue]{hyperref}
\usepackage{pdflscape}

\newtheorem{thm}{Theorem}[section]
\newtheorem{lem}[thm]{Lemma}
\newtheorem{eg}[thm]{Example}
\newtheorem{prop}[thm]{Proposition}
\newtheorem{cor}[thm]{Corollary}
\newtheorem{defn}[thm]{Definition}
\newtheorem{rem}[thm]{Remark}

\newtheorem{thm2}{Theorem}
\newtheorem{defn2}[thm2]{Definition}

\newenvironment{prf}{{\noindent \textbf{Proof:}\ }}{\hfill $\Box$\\ \smallskip}

\numberwithin{equation}{section}

\newcommand{\ti}{\tilde}

\newcommand{\smnoind}{\smallskip\noindent}

\newcommand{\Cu}{\mathcal{U}}

\newcommand{\ad}{{\rm Ad}\ \!}

\newcommand{\BC}{\mathbb{C}}

\newcommand{\BN}{\mathbb{N}}

\newcommand{\CS}{\mathcal{S}}

\newcommand{\CH}{\mathcal{H}}

\newcommand{\CL}{\mathcal{L}}

\newcommand{\CI}{\mathcal{I}}

\newcommand{\KJ}{\mathfrak{J}}
\newcommand{\KF}{\mathfrak{F}}

\newcommand{\KM}{\mathfrak{M}}
\newcommand{\KH}{\mathfrak{H}}
\newcommand{\KP}{\mathfrak{P}}

\begin{document}

\title{Spectral gap actions and invariant states}

\author{Han Li \and Chi-Keung Ng}

\address[Han Li]{Department of Mathematics, Yale University, PO Box 208283, New Haven, CT 06520-8283, USA.}
\email{li.han@yale.edu}

\address[Chi-Keung Ng]{Chern Institute of Mathematics and LPMC, Nankai University, Tianjin 300071, China.}
\email{ckng@nankai.edu.cn}

\thanks{The second named author is supported by the National Natural Science Foundation of China (11071126).}

\date{\today}

\keywords{Spectral gap, Invariant states, Property $T$, Inner amenability}

\subjclass[2010]{Primary: 46L55, 22A10, 43A65; Secondary: 22D25, 22D10, 43A07}

\begin{abstract}
We define spectral gap actions of discrete groups on von Neumann algebras and study their relations with invariant states.
We will show that a finitely generated ICC group $\Gamma$ is inner amenable if and only if there exist more than one inner invariant states on the group von Neumann algebra $L(\Gamma)$.
Moreover, a countable discrete group $\Gamma$ has property $(T)$ if and only if for any action $\alpha$ of $\Gamma$ on a von Neumann algebra $N$, every $\alpha$-invariant state on $N$ is a weak-$^*$-limit of a net of normal $\alpha$-invariant states.
\end{abstract}

\maketitle

\section{Introduction}

In this paper, $\Gamma$ is a discrete group, $K$ is a (complex) Hilbert space and $\pi$ is a unitary representation of $\Gamma$ on $K$.
We denote by $1_\Gamma$ the one dimensional trivial representation of $\Gamma$, and $K^\pi$ the set of $\pi$-invariant vectors in $K$.
Recall that (see \cite{BC}) $\pi$ is said to have a \emph{spectral gap} if the restriction of $\pi$ on the orthogonal complement $(K^\pi)^\bot$ does not weakly contain $1_\Gamma$ (in the sense of Fell), i.e.\ there does not exist a net of unit vectors $\xi_i\in (K^\pi)^\bot$ satisfying $\|\pi_t\xi_i -\xi_i\| \to 0$ for every $t\in \Gamma$.

The following result on group actions, which relates the existence of spectral gap of certain representation and the uniqueness of the invariant states (means), is crucial to the solution of the Banach-Ruziewicz problem (see \cite{Marg},\cite{Sull},\cite{Drin}).

\begin{thm2}\label{motivation}
Let $\Gamma$ be a countable group with a measure preserving ergodic action on a standard non-atomic probability space $(X, \mu)$, and $\rho$ be the associated unitary representation of $\Gamma$ on $L^2(X, \mu)$.

\smnoind
(a) (\cite{dR}, \cite{Sch}, \cite{LH}) If $\rho$ has a spectral gap, then the integration with respect to $\mu$ is the unique $\Gamma$-invariant state on $L^\infty(X,\mu)$.

\smnoind
(b) (\cite{Sch}) If there is a unique $\Gamma$-invariant state on $L^\infty(X,\mu)$, then $\rho$ has a spectral gap.
\end{thm2}

Motivated by this result we will define spectral gap actions of discrete groups on von Neumann algebras, and study the corresponding invariant states on the von Neumann algebras.

From now on, $N$ is a von Neumann algebra  with standard form $(N, \KH, \KJ, \KP)$, i.e.\ $\KH$ is a Hilbert space, $\KJ:\KH \to \KH$ is a conjugate linear bijective isometry, and $\KP$ is a self-dual cone in $\KH$ such that there is a (forgettable) faithful representation of $N$ on $\KH$ satisfying some compatibility conditions with $\KJ$ and $\KP$.
The readers may find in \cite{Haag-st-form} more information on this topic.
As proved in \cite[Theorem 3.2]{Haag-st-form}, there exists a unique unitary representation $\Cu$ of the $^*$-automorphism group ${\rm Aut}(N)$ on $\KH$ such that
$$g(x) = \Cu_{g} x \Cu_{g}^*, \quad \Cu_g\KJ=\KJ\Cu_g \quad \text{and} \quad \Cu_g(\KP)\subseteq \KP\qquad (x\in N; g\in{\rm Aut}(N)).$$
Hence, an action $\alpha$ of $\Gamma$ on $N$ gives rise to a unitary representation $\Cu_{\alpha}=\Cu\circ\alpha$ of $\Gamma$ on $\KH$.

\begin{defn2}
Let $\alpha$ be an action of $\Gamma$ on $N$.

\smnoind
(a) $\alpha$ is said to have a \emph{spectral gap} if the representation $\Cu_{\alpha}$ has a spectral gap.

\smnoind
(b) $\alpha$ is said to be \emph{standard} if every $\alpha$-invariant state on $N$ is a weak-$^*$-limit of a net of normal $\alpha$-invariant states.
\end{defn2}

We will study the relationship between spectral gap actions and standard actions of discrete groups on von Neumann algebras. This is the main topic of Section 2. Our first main result generalizes Theorem \ref{motivation}(a) to the situation of discrete group actions on von Neumann algebras.

\begin{thm2}\label{thm:main}Let $\alpha$ be an action of a (possibly uncountable) discrete group $\Gamma$ on a von Neumann algebra $N$. If $\alpha$ has a spectral gap, then $\alpha$ is standard.
\end{thm2}

In contrast to Theorem \ref{motivation}(b), the converse statement of Theorem \ref{thm:main} is not true in general (see Example \ref{eg:4-not-2}(b)).
Nevertheless, we will give several situations in which the converse  does hold (see e.g.\ Proposition \ref{prop:inv-mean-sp-gap-act}(b)).
Let us also note that Theorem \ref{motivation} (and hence Theorem \ref{thm:main}) does not extend to locally compact groups.
Indeed, if $G$ is the circle group (equipped with the canonical compact topology) and $N := L^\infty(G)$ with the left translation action $\alpha$ by $G$, then there is a unique normal $\alpha$-invariant state on $N$, but there exist more than one $\alpha$-invariant states (see e.g.\ \cite[Proposition 2.2.11]{Lub}).

We will give two applications of Theorem \ref{thm:main}.
The first one concerns inner amenability and will be considered in Section 3.
The notion of inner amenability was first introduced by
Effros, aiming to give a group theoretic description of property Gamma of the group von Neumann algebra $L(\Gamma)$ of an ICC group.
Recall that  that
$\Gamma$ is an \emph{ICC group} if $\Gamma \neq \{e\}$ and all the non-trivial conjugacy classes of $\Gamma$ are infinite.
Moreover, $\Gamma$ is \emph{inner amenable} if there exist more than one inner invariant states on $\ell^\infty(\Gamma)$.
Effros showed in \cite{Eff} that a countable ICC group $\Gamma$ is inner amenable if $L(\Gamma)$ has property Gamma.
However, Vaes recently gave, in \cite{Vaes}, a counter example to the converse.
From an opposite angle, it is natural to ask whether one can express inner amenability of $\Gamma$ in terms of certain property of $L(\Gamma)$.
One application of Theorem \ref{thm:main} is  the following result.

\begin{thm2}\label{thm:inner}
Let $\Gamma$ be a finitely generated ICC group. Then $\Gamma$ is inner amenable if and only if there exist more than one inner invariant states on the group von Neumann algebra $L(\Gamma)$.
\end{thm2}

We will also study an alternative generalization of inner amenability (which is called ``strongly inner amenability'') of  (not necessarily ICC) discrete groups, which is of independent interest.

Section 4 is concerned with our second application, which is related to property $(T)$.
Recall that $\Gamma$ is said to have \emph{property $(T)$} if every unitary representation of $\Gamma$ has a spectral gap.
It follows from Theorem \ref{thm:main} that if $\Gamma$ has property $(T)$, then every action of $\Gamma$ on a von Neumann algebra is standard, and in particular, the absence of normal invariant state implies the absence of invariant state for this action.
We will show that this property actually characterizes property $(T)$.

\begin{thm2}\label{thm:property-T} Let $\Gamma$ be a countable discrete group.
The following statements are equivalent.
\begin{enumerate}
\item $\Gamma$ has property $(T)$.
\item All actions of $\Gamma$ on von Neumann algebras are standard.
\item For any action $\alpha$ of $\Gamma$ on a von Neumann algebra with only one normal $\alpha$-invariant state, there is only one $\alpha$-invariant state.
\item For every action $\alpha$ of $\Gamma$ on a von Neumann algebra without normal $\alpha$-invariant state, there is no $\alpha$-invariant state.
\end{enumerate}
\end{thm2}

We will also consider the implication $(4)\Rightarrow (1)$ for general discrete group  $\Gamma$ with property $(T, FD)$ in the sense of \cite{LZ89} (Proposition \ref{prop:no-inv-state}). Consequently, a minimally almost period discrete group satisfying (4) is finitely generated.

We end this introduction by introducing more notation.
Throughout this article, $\overline{K}$ and $\CL(K)$ denote the conjugate Hilbert space of $K$ and the space of bounded linear maps on $K$ respectively.
For any $\xi,\eta\in K$, we define $\omega_{\xi,\eta}\in \CL(K)^*$ by $\omega_{\xi,\eta}(x) := \langle x\xi, \eta\rangle$ ($x\in \CL(K)$) and set $\omega_\xi := \omega_{\xi,\xi}$.
We consider $p_\pi\in \CL(K)$ to be the orthogonal projection onto $K^\pi$.
Furthermore, we denote by $C^*(\pi)$ the $C^*$-subalgebra generated by $\{\pi_t\}_{t\in \Gamma}$ and by $vN(\pi)$ the bicommutant of $C^*(\pi)$.
We define a functional $\varepsilon_\pi: C^*(\pi) \to \BC$ formally by $\varepsilon_\pi\big(\sum_{k=1}^n c_k \pi_{t_k}\big) = \sum_{k=1}^n c_k$ (warning: $\varepsilon_\pi$ is not necessarily well-defined.)

\section{Invariant States and Spectral Gap Actions (Theorem \ref{thm:main})}

Let us start with the following lemma, which may be a known result.
Note that part (d) of this lemma is motivated by \cite{Pasch} (and one can also obtain this part using the argument in \cite{Pasch}).

\begin{lem}\label{lem:p_pi}
(a) $\varepsilon_\pi$ is well-defined if and only if $1_\Gamma$ is weakly contained in $\pi$.
This is equivalent to the existence of $\psi\in \CL(K)^*_+$ satisfying
\begin{equation}\label{eqt:m-pi}
\psi(\pi_t)\ =\ 1 \qquad (t\in \Gamma)
\end{equation}
(in this case, $\varepsilon_\pi = \psi|_{C^*(\pi)}$).

\smnoind
(b) $\pi$ does not have a spectral gap if and only if there is $\psi\in \CL(K)_+^*$ satisfying \eqref{eqt:m-pi} and $\psi(p_{\pi}) = 0$.

\smnoind
(c) If $\psi\in \CL(K)_+^*$ satisfies \eqref{eqt:m-pi}, then $\psi(\pi_t x) = \psi(x)$ ($t\in \Gamma, x\in \CL(K)$) and $\psi$ is $\ad\pi$-invariant.

\smnoind
(d) $p_\pi\notin C^*(\pi)$ if and only if $p_\pi \neq 0$ and $\pi$ does not have a spectral gap.
\end{lem}
\begin{prf}
Let $\pi^0$ be the restriction of $\pi$ on $(K^\pi)^\bot$.
One has $\pi_t = 1_{K^\pi}\oplus \pi^0_t$ ($t\in \Gamma$) and $p_\pi = 1_{K^\pi}\oplus 0$.

\smnoind
(a) For any unitary representation $\rho$ of $\Gamma$, we denote by $\ti \rho$ the induced $^*$-representation of the full group $C^*$-algebra of $\Gamma$.
Then $\varepsilon_\pi$ is well-defined if and only if $\ker \ti\pi \subseteq \ker \ti 1_\Gamma$, which in turn is equivalent to $1_\Gamma$ being weakly contained in $\pi$ (see e.g.\ \cite[Theorem F.4.4]{BHV}).
The second statement is trivial.

\smnoind
(b) Part (a) implies that $\pi$ does not have a spectral gap if and only if there is $\phi\in \CL\big((K^\pi)^\bot\big)^*_+$ satisfying $\phi(\pi_t^0) = 1$ for each $t\in \Gamma$.
This, in turn, is equivalent to the existence of $\varphi\in \big(\CL(K^\pi)\oplus \CL((K^\pi)^\bot)\big)^*_+$ with $\varphi(\pi_t) = 1$ ($t\in \Gamma$) and $\varphi(p_\pi) = 0$.

\smnoind
(c) This can be obtained by applying the Cauchy-Schwarz inequality to $\psi\big((1-\pi_t)x\big)$.

\smnoind
(d) Notice that any of the two statements implies that $p_\pi \neq 0$, and so $\varepsilon_\pi$ is well-defined (because of part (a)).

\smnoind
$\Rightarrow)$.
Since $p_\pi\notin C^*(\pi)$, the map $\rho:C^*(\pi) \to \CL(K)$ defined by  $\rho(x) := x - \varepsilon_\pi(x) p_\pi$ is injective with its image being $0\oplus C^*(\pi^0)$.
Consequently, $\varepsilon_{\pi^0} = \varepsilon_\pi\circ \rho^{-1}$ (we identify $0\oplus C^*(\pi^0)$ with $C^*(\pi^0)$) and is well-defined.
This shows that $1_\Gamma$ is weakly contained in $\pi^0$ (by part (a)).

\smnoind
$\Leftarrow)$.
By parts (a) and (b), there is $\psi\in \CL(K)^*_+$ with $\psi|_{C^*(\pi)} = \varepsilon_\pi$ and $\psi(p_\pi) = 0$.
Take any unit vector $\xi\in K^\pi$.
If $p_\pi \in C^*(\pi)$, then $1 = \omega_\xi(p_\pi) = \varepsilon_\pi(p_\pi) = 0$ (as $\varepsilon_\pi = \omega_\xi$ on $C^*(\pi)$), which is absurd.
\end{prf}

Let $N$ be a von Neumann algebra.
We denote by $N_*$ the pre-dual space of $N$, which is naturally identified with the subspace of all weak-$^*$-continuous linear functionals in the dual space $N^*$.
Let
$$N^*_\alpha:= \{f\in N^*: f\circ \alpha_t =f \text{ for any } t\in \Gamma\}, \quad \KM_{N,\alpha} := \{f\in N^*_\alpha: f\geq 0 \text{ and } f(1) =1\}$$
and $\KM^{N,\alpha} := \overline{\KM_{N,\alpha}\cap N_*}^{\sigma(N^*,N)}$ (where $\overline{E}^{\sigma(N^*,N)}$ means the weak-$^*$-closure of a subset $E\subseteq N^*$).
Clearly, $\KM^{N,\alpha}\subseteq \KM_{N,\alpha}$.

The following theorem is a refined version of Theorem \ref{thm:main}.
This result can also be regarded as generalizations of the main results in \cite{Pasch} and in \cite{LH} (see also \cite[Proposition 3.4.1]{Lub}).

\begin{thm}\label{thm:sp-gap-act-inv-mean}
Let $\alpha$ be an action of a group $\Gamma$ on a von Neumann algebra $N$  with standard form $(N, \KH, \KJ, \KP)$.
Consider the following statements.

\begin{enumerate}[(G1)]
\item If $\psi\in \CL(\KH)_+^*$ satisfies $\psi(\Cu_{\alpha_t}) = 1$ ($t\in \Gamma$), then $\psi(p_{\Cu_{\alpha}}) \neq 0$.
\item $\alpha$ has a spectral gap.
\item $p_{\Cu_{\alpha}}\in C^*(\Cu_{\alpha})$
\item $\KM_{N,\alpha} = \KM^{N,\alpha}$.
\end{enumerate}

\noindent
One has $(G1) \Leftrightarrow (G2) \Rightarrow (G3)$ and $(G2)\Rightarrow (G4)$.
\end{thm}
\begin{prf}
$(G1) \Leftrightarrow (G2)$.
This follows from Lemma \ref{lem:p_pi}(b).

\smnoind
$(G2) \Rightarrow (G3)$.
This follows from Lemma \ref{lem:p_pi}(d).

\smnoind
$(G2) \Rightarrow (G4)$.
Let $m\in N^*_+$ be an $\alpha$-invariant state.
If $\{f_i\}_{i\in I}$ is a net of states in $N_*$ that $\sigma(N^*, N)$-converges to $m$, the ``convergence to invariance'' type argument (see e.g.\ \cite[p.33-34]{Green}) will produce a net $\{g_j\}_{j\in J}$ in the convex hull of $\{f_i\}_{i\in I}$ that $\sigma(N^*, N)$-converges to $m$ and satisfies $\|\alpha_t^*(g_j) - g_j\|_{N_*} \to 0$ for every $t\in \Gamma$.

For each $j\in J$, there is a unique unit vector $\zeta_j\in \KP$ with $g_j = \omega_{\zeta_j}$.
As $\alpha_{t^{-1}}^*(g_j) = \omega_{\Cu_{\alpha_{t}}(\zeta_j)}$,
we have, by \cite[Lemma 2.10]{Haag-st-form} (and the fact that $\Cu_{\alpha_{t}}(\KP)\subseteq \KP$),
\begin{equation}\label{eqt:Pow-Stor-ineq}
\|\Cu_{\alpha_{t}}(\zeta_j) - \zeta_j\|^2
\ \leq\  \|\alpha_{t^{-1}}^*(g_j) - g_j\|_{N_*}
\qquad (t\in \Gamma).
\end{equation}

Let $\zeta^0_j := \zeta_j - p_{\Cu_{\alpha}}(\zeta_j)$.
If $\|\zeta_j^0\| \nrightarrow 0$, a subnet of $\{\zeta_j^0\}_{j\in J}$ will produce an almost $\Cu_{\alpha}$-invariant unit vector in  $(\KH^{\Cu_{\alpha}})^\bot$, which contradicts the spectral gap assumption of $\Cu_{\alpha}$.
Consequently, if we set $\zeta_j^1:= \frac{p_{\Cu_{\alpha}}(\zeta_j)}{\|p_{\Cu_{\alpha}}(\zeta_j)\|}$, then $\{\omega_{\zeta_j^1}\}_{j\in J}$ is a net of $\alpha$-invariant states that  $\sigma(N^*, N)$-converges to $m$.
\end{prf}

It is easy to see that $\KM_{N,\alpha}$ spans $N^*_\alpha$ (see e.g.\ \cite[Proposition 2.2]{Pat}).
Moreover, if we denote $N_*^\alpha :=N^*_\alpha\cap N_*$, it is not hard to check that $\KM^{N,\alpha}$ spans $\overline{N_*^\alpha}^{\sigma(N^*, N)}$.
Thus, Statement (G4) is equivalent to the fact that $N^*_\alpha = \overline{N_*^\alpha}^{\sigma(N^*, N)}$.

\begin{lem}\label{lem:H-N_*}
(a) The map sending $\xi \in \KH$ to $\omega_\xi \in N_*$ restricts to a bijection from $\KP \cap \KH^{\Cu_{\alpha}}$ onto the positive part of
$N_*^\alpha$.

\smnoind
(b) If $\KH^{\Cu_{\alpha}}$ is finite dimensional, then $N_*^\alpha$ is finite dimensional.

\smnoind
(c) If $\dim \KH^{\Cu_{\alpha}} =1$, then $\dim N_*^\alpha = 1$.
\end{lem}
\begin{prf}
(a) This follows from the fact that for any $g\in (N_*^\alpha)_+$,  if $\zeta\in \KP$ is the unique element satisfying $g = \omega_\zeta$, we have $\zeta\in \KH^{\Cu_{\alpha}}$, by the inequality in \eqref{eqt:Pow-Stor-ineq}.

\smnoind
(b) For any $\phi\in N_*^\alpha$, there exist $\alpha$-invariant elements $\phi_1,\phi_2,\phi_3,\phi_4\in (N_*^\alpha)_+$ with $\phi = (\phi_1 - \phi_2) + \mathrm{i} (\phi_3 - \phi_4)$.
By part (a), one can find $\zeta_k\in \KP\cap \KH^{\Cu_{\alpha}}$ such that $\phi_k = \omega_{\zeta_k}$ ($k=1,2,3,4$).
Therefore, the linear map from $(\KH^{\Cu_{\alpha}} \otimes \overline{\KH^{\Cu_{\alpha}}}) \oplus (\KH^{\Cu_{\alpha}} \otimes \overline{\KH^{\Cu_{\alpha}}})$ to $N_*^\alpha$ given by $(\xi_1\otimes \overline{\eta_1}, \xi_2\otimes \overline{\eta_2}) \mapsto \omega_{\xi_1,\eta_1} + \mathrm{i} \omega_{\xi_2,\eta_2}$ is surjective and $N_*^\alpha$ is finite dimensional.

\smnoind
(c) Clearly, $N_*^\alpha \neq (0)$.
If $\dim N_*^\alpha > 1$, there will be two different norm one elements in $(N_*^\alpha)_+$.
Thus, part (a) gives two norm one elements in $\KP \cap \KH^{\Cu_{\alpha}}$, which is not possible.
\end{prf}

The following proposition concerns some converses of the implications in Theorem \ref{thm:sp-gap-act-inv-mean}.
Note that the idea of the argument for $p_{\Cu_{\alpha}}\in vN(\Cu_{\alpha})$ in part (b) comes from \cite{Pasch}.

\begin{prop}\label{prop:inv-mean-sp-gap-act}
Let $\alpha$ be an action of a group $\Gamma$ on a von Neumann algebra $N$.

\smnoind
(a) If $p_{\Cu_{\alpha}}\neq 0$, then Statement (G3) implies Statement (G1).

\smnoind
(b) If $p_{\Cu_{\alpha}}\in N$, then Statement (G4) implies Statement (G1) and $p_{\Cu_{\alpha}}$ belongs to the center $Z\big(vN(\Cu_{\alpha})\big)$ of $vN(\Cu_{\alpha})$.
\end{prop}
\begin{prf}
(a) This follow from Lemma \ref{lem:p_pi}(d).

\smnoind
(b) Suppose that Statement (G4) holds, but there exists $\psi\in \CL(\KH)^*_+$ with $\psi(\Cu_{\alpha_t}) = 1$ ($t\in \Gamma$) and $\psi(p_{\Cu_{\alpha}}) = 0$.
Then Lemma \ref{lem:p_pi}(c) implies that $\psi|_{N}$ is $\alpha$-invariant, and Statement (G4) gives a net $\{g_j\}_{j\in J}$ of states in $N_*^\alpha$ that $\sigma(N^*,N)$-converges to $\psi|_N$.
By Lemma \ref{lem:H-N_*}(a), one can find $\zeta_j\in \KP\cap \KH^{\Cu_{\alpha}}$ with $g_j = \omega_{\zeta_j}$, which means that $g_j(p_{\Cu_{\alpha}}) = 1$ ($j\in J$).
This contradicts $\psi(p_{\Cu_{\alpha}}) = 0$.

For the second conclusion, it is clear that $p_{\Cu_{\alpha}}\in vN(\Cu_{\alpha})'$ because $\KH^{\Cu_{\alpha}}$ is $\Cu_{\alpha}$-invariant.
Suppose on the contrary that $p_{\Cu_{\alpha}}\notin vN(\Cu_{\alpha})$.
Since $\BC p_{\Cu_{\alpha}}$ is one dimensional and
\begin{equation}\label{eqt:rel-p}
\Cu_{\alpha_t} p_{\Cu_{\alpha}}\ =\ p_{\Cu_{\alpha}} \qquad (t\in \Gamma),
\end{equation}
one sees that $vN(\Cu_{\alpha}) + \BC p_{\Cu_{\alpha}}$ is a von Neumann algebra.
Moreover, if $\xi\in \KH^{\Cu_{\alpha}}$ is a norm one vector (note that we assume $p_{\Cu_{\alpha}}\notin vN(\Cu_{\alpha})$), we define a functional on $vN(\Cu_{\alpha}) + \BC p_{\Cu_{\alpha}}$ by
$$\phi(x+ c p_{\Cu_{\alpha}})\ :=\ \omega_{\xi}(x)
\qquad (x\in vN(\Cu_{\alpha}); c\in \BC).$$
Clearly, $\phi$ is weak-$^*$-continuous and it is a $^*$-homomorphism since $xp_{\Cu_{\alpha}}\in \BC p_{\Cu_{\alpha}}$ (because of \eqref{eqt:rel-p}).
If $\psi$ is a normal state extension of $\phi$ on $\CL(\KH)$, then $\psi(p_{\Cu_{\alpha}}) = 0$ and Lemma \ref{lem:p_pi}(c) implies that $\psi|_N$ is an $\alpha$-invariant normal state.
Now, Lemma \ref{lem:H-N_*}(a) produces $\zeta\in \KP\cap \KH^{\Cu_{\alpha}}$ with $\psi|_N = \omega_{\zeta}$ and we have the contradiction that $\psi(p_{\Cu_{\alpha}}) = 1$.
\end{prf}

Our next example shows that Statement (G4) does not imply Statement (G2) in general.
We first set some notation and recall some facts.
Let $\CI_N$ be the inner automorphism group of $N$ and $\beta$ be the canonical action of $\CI_N$ on $N$.
Then $\beta$-invariant states on $N$ are precisely tracial states.
Suppose, in addition, that $N$ admits a normal faithful tracial state $\tau$.
If $(H_\tau, \Psi_\tau)$ is the GNS construction with respect to $\tau$ and $\Lambda_\tau: N \to H_\tau$ is the canonical map, then $\KH \cong H_\tau$ and the canonical action of $N$ on $\KH$ can be identified with $\Psi_\tau$. 
For any $g\in {\rm Aut}(N)$, one has $\Cu_{\beta_g}(\Lambda_\tau (x)) = \Lambda_\tau(g(x))$ ($x\in N$).

\begin{eg}\label{eg:4-not-2}
Suppose that $\lambda$ is the left regular representation of $\Gamma$ and put $L(\Gamma):= vN(\lambda)$.
Then $\KH = \ell^2(\Gamma)$ and we let $\CI_{L(\Gamma)}$ and $\beta$ be as in the above.

\smnoind
(a) The representation $\Cu_{\beta}\circ {\rm Ad}\circ \lambda: \Gamma\to \CL(\ell^2(\Gamma))$ coincides with the ``conjugate representation'' $\gamma$ defined by  $\gamma_t(\xi)(s) := \xi (t^{-1}st)$ ($s,t\in \Gamma; \xi\in \ell^2(\Gamma)$).

\smnoind
(b) Suppose that $\Gamma$ is an amenable countable ICC group.
Since $L(\Gamma)$ is a type $I\!I_1$-factor, it has only one tracial state and this state is normal.
Consequently, $\KM_{L(\Gamma), \beta} = \KM^{L(\Gamma), \beta}$ (because $\beta$-invariant states are precisely tracial states).

On the other hand, we have $\KH^{\Cu_\beta} = \BC\delta_e$. 
Moreover, as $L(\Gamma)$ is semi-discrete, it has property Gamma (see \cite[Corollary 2.2]{Con}), and the restriction of $\Cu_{\beta}$ on $(\KH^{\Cu_\beta})^\bot$ weakly contains the trivial representation (by an equivalent form of property Gamma in \cite[Theorem 2.1(c)]{Con}).
Thus, the representation $\Cu_\beta$ (and hence the action $\beta$) does not have a spectral gap.
\end{eg}

\section{Invariant States and Inner Amenability (Theorem \ref{thm:inner})}

In this section, we will give an application of Theorem \ref{thm:sp-gap-act-inv-mean} to inner amenability.
We recall (from the main theorem in \cite{Eff}) that an ICC group $\Gamma$ is inner amenable if and only if there is a net $\{\xi_i\}_{i\in I}$ of unit vectors in $\ell^2(\Gamma\setminus \{e\})$ such that $\|\gamma_t(\xi_i) - \xi_i\| \to 0$ for any $t\in \Gamma$ (where $\gamma$ is as in Example \ref{eg:4-not-2}(a)).
Notice that this inner amenability is slightly different from the one in \cite{LP} and \cite{Pat}, and it is called ``non-trivially inner amenable'' in \cite[p.84]{Pat}.

Let us consider another extension of inner amenability to general (not necessarily ICC) discrete groups.

\begin{defn}
$\Gamma$ is said to be \emph{strongly inner amenable} if the conjugate representation $\gamma$ does not have a spectral gap.
\end{defn}

It is obvious that strong inner amenability implies inner amenability, and the converse holds when $\Gamma$ is an ICC group.
On the other hand, all abelian groups and all property $(T)$ groups are not strongly inner amenable.

\begin{prop}\label{prop:st-inner}
Let $\alpha$ be the action of $\Gamma$ on $\CL(\ell^2(\Gamma))$ given by $\alpha_t(x) := \gamma_t x \gamma_{t^{-1}}$, $\Gamma_{\rm fin}$ be the normal subgroup consisting of elements in $\Gamma$ with finite conjugacy classes and $A(\Gamma)$ be the predual of $L(\Gamma)$.

\smnoind
(a) Consider the following statements.

\begin{enumerate}[(S1)]
\item $\Gamma/\Gamma_{\rm fin}$ is inner amenable.
\item $\ell^\infty(\Gamma)^*_\alpha$ is not the $\sigma(\ell^\infty(\Gamma)^*, \ell^\infty(\Gamma))$-closure of $\ell^1(\Gamma)^\alpha$.
\item $\Gamma$ is strongly inner amenable.
\item There is $\psi\in \CL(\ell^2(\Gamma))_+^*$ satisfying $\psi(\gamma_t) = 1$ ($t\in \Gamma$) and $\psi(p_\gamma) = 0$.
\item $p_\gamma\notin C^*(\gamma)$.
\item $L(\Gamma)^*_\alpha$ is not the $\sigma(L(\Gamma)^*, L(\Gamma))$-closure of $A(\Gamma)^\alpha$.
\item $\KM_{L(\Gamma),\alpha}$ does not coincide with the set of tracial states on $L(\Gamma)$.
\end{enumerate}

\noindent
Then $(S1) \Rightarrow (S2) \Rightarrow (S3) \Leftrightarrow (S4) \Leftrightarrow (S5)$ and  $(S6)\Leftrightarrow (S7) \Rightarrow (S3)$.

\smnoind
(b) If $\Gamma_{\rm fin}$ coincides with the center of $\Gamma$, then $(S5) \Rightarrow (S1)$ and $p_\gamma\in Z(vN(\gamma))$.

\smnoind
(c) If $\Gamma_{\rm fin}$ is finite, we have $(S5) \Rightarrow (S1)$.
\end{prop}
\begin{prf}
Notice that if $x\in L(\Gamma)$, then $\alpha_t(x) = \lambda_t x\lambda_t^*\in L(\Gamma)$.
On the other hand, if $y\in \ell^\infty(\Gamma)$ (considered as a subalgebra of $\CL(\ell^2(\Gamma))$ in the canonical way), then $\alpha_t(y)\in \ell^\infty(\Gamma)$ and $\alpha_t(y)(s) = y(t^{-1}st)$.
Moreover, if we either set
$$N = \ell^\infty(\Gamma) \quad \text{or} \quad N = L(\Gamma)$$
(and we denote by $\alpha$ the restriction of $\alpha$ on $N$ by abuse of notation), then $\KH = \ell^2(\Gamma)$ and $\Cu_{\alpha}$ coincides with $\gamma$ (even though $\KJ$ and $\KP$ are different in these two cases).

\smnoind
(a) By Theorem \ref{thm:sp-gap-act-inv-mean} and Proposition \ref{prop:inv-mean-sp-gap-act}(a) (for both $N=\ell^\infty(\Gamma)$ and $N= L(\Gamma)$), it remains to show that $(S1)\Rightarrow (S2)$ and $(S6)\Leftrightarrow (S7)$.

\noindent
$(S1)\Rightarrow (S2)$.
Since $\Gamma/\Gamma_{\rm fin}$ is inner amenable, we know from \cite[Corollary 1.4]{Kan-Mark} as well as the argument of Lemma \ref{lem:p_pi}(c) that there is an $\alpha$-invariant mean $m$ on $\Gamma$ such that $m(\chi_{\Gamma_{\rm fin}}) \neq 1$.
Obviously, $\ell^2(\Gamma)^\gamma \subseteq \ell^2(\Gamma_{\rm fin})$ and Lemma \ref{lem:H-N_*}(a) implies that $g(\chi_{\Gamma_{\rm fin}}) = \|g\|$ if $g\in \ell^1(\Gamma)^\alpha_+$.
Consequently, $m\notin \KM^{\ell^\infty(\Gamma),\alpha}$.

\noindent
$(S6)\Rightarrow (S7)$.
Assume that $\KM_{L(\Gamma),\alpha}$ coincides with the set of all tracial states on $L(\Gamma)$.
It is well-known that every tracial state on $L(\Gamma)$ is a weak-$^*$-limit of a net of normal tracial states (see e.g.\ \cite[Proposition 8.3.10]{KRII}).
Thus, $\KM_{L(\Gamma),\alpha} = \KM^{L(\Gamma),\alpha}$, which contradicts Statement (S6).

\noindent
$(S7)\Rightarrow (S6)$.
Clearly, any tracial state on $L(\Gamma)$ is $\alpha$-invariant.
Moreover, $\omega\in A(\Gamma)$ is $\alpha$-invariant if and only if $$\omega(\lambda_t x)\ =\ \omega(x\lambda_{t}) \qquad (x\in L(\Gamma); t\in \Gamma),$$
which is equivalent to $\omega(yx) = \omega(xy)$ ($x,y\in L(\Gamma)$).
Hence, $\KM_{L(\Gamma),\alpha}\cap A(\Gamma)$ is the set of all normal tracial states on $L(\Gamma)$.
Consequently, if $L(\Gamma)^*_\alpha$ is the weak-$^*$-closure of $A(\Gamma)^\alpha$, then $\KM_{L(\Gamma),\alpha}$ coincides with the set of tracial states, which contradicts Statement (S7).

\smnoind
(b) The assumption implies that every finite conjugacy class of $\Gamma$ is a singleton set.
Consequently, one has $\ell^2(\Gamma)^\gamma = \ell^2(\Gamma_{\rm fin})$ and $p_\gamma = \chi_{\Gamma_{\rm fin}}\in \ell^\infty(\Gamma)$.
Thus, $p_\gamma\in Z(vN(\gamma))$ by Proposition \ref{prop:inv-mean-sp-gap-act}(b).
Furthermore, Statement (S4) produces an $\alpha$-invariant mean $m$ on $\Gamma$ (see Lemma \ref{lem:p_pi}(c)) such that $m(p_\gamma) =0$.
Now, by \cite[Corollary 1.4]{Kan-Mark}, $\Gamma/\Gamma_{\rm fin}$ is inner amenable.

\smnoind
(c) Let $C_1,...,C_n$ be all the finite conjugacy classes of $\Gamma$.
Let $p_k\in \CL(\ell^2(\Gamma))$ be the orthogonal projection onto $\ell^2(C_k)$ ($k=1,...,n$) and $p_0\in \CL(\ell^2(\Gamma))$ be the orthogonal projection onto $\ell^2(\Gamma\setminus \Gamma_{\rm fin})$.
Then
$$\ell^2(\Gamma)^\gamma = \{\xi\in \ell^2(\Gamma): p_0\xi = 0 \text{ and } p_k\xi \text{ is a constant function for every } k=1,...,n\}.$$
Suppose that $\{\xi_i\}_{i\in I}$ is an almost $\gamma$-invariant unit vector in $(\ell^2(\Gamma)^\gamma)^\bot$.
One can find a subnet $\{\xi_{i_j}\}_{j\in J}$ such that $\{p_k\xi_{i_j}\}_{j\in J}$ is norm-converging to some $\xi^{(k)} \in \ell^2(C_k)$ with $\sum_{s\in C_k}\xi^{(k)}(s) = 0$ (for $k =1,...,n$).
As $\|\gamma_t(p_k \xi_{i_j}) - p_k \xi_{i_j}\|\to 0$, we see that $\gamma_t \xi^{(k)} = \xi^{(k)}$ ($t\in \Gamma$), which means that $\xi^{(k)} = 0$.
Consequently, $\|p_0\xi_{i_j}\| \to 1$ and a subnet of $\{\omega_{p_0\xi_{i_j}}\}_{j\in J}$ will produce a $\gamma$-invariant mean $m$ satisfying $m(\chi_{\Gamma_{\rm fin}}) =0$.
Now, Statement (S1) follows from \cite[Corollary 1.4]{Kan-Mark}.
\end{prf}

The following theorem concerns the case when there is a unique inner invariant state on $L(\Gamma)$.
Part (b) of which gives Theorem \ref{thm:inner}.

\begin{thm}\label{thm:ICC}
Let $\alpha$ be as in Proposition \ref{prop:st-inner}.

\smnoind
(a) If $\Gamma$ is amenable, then there are more than one $\alpha$-invariant states on $L(\Gamma)$.

\smnoind
(b) Suppose that $\Gamma$ is an ICC group.
If $\Gamma$ is not inner amenable, there is a unique $\alpha$-invariant state on $L(\Gamma)$.
The converse holds when $\Gamma$ is finitely generated.
\end{thm}
\begin{prf}
(a) Suppose that $n\in \ell^\infty(\Gamma)^*_+$ is an invariant mean.
Using the ``convergence to invariance'' argument, there is a net $\{g_j\}_{j\in J}$ of norm one elements in $\ell^1(\Gamma)_+$ such that $\|t\cdot g_j - g_j\|_{\ell^1(\Gamma)} \to 0$ and $\|g_j\cdot t- g_j\|_{\ell^1(\Gamma)} \to 0$, for every $t\in \Gamma$ (where $(t\cdot g_i)(s) = g_i(t^{-1}s)$ and  $(g_i\cdot t)(s) = g_i(st)$).
Set $\eta_j := g_j^{1/2}\in \ell^2(\Gamma)$ ($j\in J$).
Then
\begin{equation*}\label{eqt:left-and-adj-inv}
\|\lambda_t(\eta_j) - \eta_j\|_{\ell^2(\Gamma)} \to 0
\qquad \text{and}\qquad
\|\gamma_t(\eta_j) - \eta_j\|_{\ell^2(\Gamma)} \to 0,
\end{equation*}
If $\phi$ is the $\sigma(L(\Gamma)^*, L(\Gamma))$-limit of a subnet of $\{\omega_{\eta_j}\}_{j\in J}$, then $\phi$ is an $\alpha$-invariant state on $L(\Gamma)$ but $\phi\neq \omega_{\delta_e}$ because $\phi(\lambda_t) = 1$ for all $t\in \Gamma$.

\smnoind
(b) Since $\Gamma$ is an ICC group, one has $\KH^{\Cu_{\alpha}} = \BC \delta_e$ and $A(\Gamma)^\alpha = \BC \omega_{\delta_e}$.
The first statement follows from Proposition \ref{prop:st-inner}(a).

To show the second statement, we suppose, on the contrary, that $\Gamma$ has a finite generating subset $\{t_1,...,t_n\}$, $\omega_{\delta_e}$ is the only $\alpha$-invariant state on $N:=L(\Gamma)$, but $\gamma$ does not have a spectral gap.
Then there is a sequence $\{\xi_k\}_{k\in \BN}$ of unit vectors in $\ell^2(\Gamma)$ satisfying
$$\langle \xi_k, \delta_e\rangle\ =\ 0 \qquad \text{and} \qquad  \|\lambda_{t_j}(\rho_{t_j}(\xi_k)) -\xi_k\| \to 0 \quad (j=1,...,n)$$
(notice that $\gamma = \lambda\circ \rho$, where $\rho$ is the right regular representation of $\Gamma$).
Set $u_j := \lambda_{t_j}\in N$ ($j=1,...,n$) and consider $\KF$ to be a free ultrafilter on $\BN$.
Since the relative commutant of $F:=\{u_1,...,u_n\}$ in $N$ is the same as the relative commutant of $F$ in the ultrapower $N^\KF$ (by the argument of \cite[Lemma 2.6]{Con}; see lines 3 and 4 of \cite[p.86]{Con}), the relative commutant of $F$ in $N^\KF$ is $\BC$.
Now, one may employ the argument for Case (1) of ``(c) $\Rightarrow$ (b)'' in \cite[p.87]{Con} to conclude that there exists a subsequence of $\{\omega_{\xi_k}\}_{k\in \BN}$ which is not relatively $\sigma(N_*,N)$-compact in $N_*$.
Hence, this subsequence has a subnet that $\sigma(N^*,N)$-converges to a non-normal state $\varphi\in N^*$.
As $\varphi$ is invariant under $\{\ad \gamma_{t_1}, ..., \ad \gamma_{t_n}\}$, it is $\alpha$-invariant and we obtain a contradiction.
\end{prf}

\begin{rem}\label{rem:icc}
(a) It was shown in \cite{Pasch} that inner amenability of an ICC group $\Gamma$ is equivalent to $p_\gamma\notin C^*(\gamma)$.
However, Proposition \ref{prop:st-inner}(b) tells us that $p_\gamma\in vN(\gamma)$ for any ICC group $\Gamma$, whether or not it is inner amenable.

\smnoind
(b) A functional $\varphi\in L(\Gamma)^*$ is $\alpha$-invariant if and only if $\varphi(ax) = \varphi(xa)$ ($x\in L(\Gamma);a\in C^*(\lambda) = C^*_r(\Gamma)$).
Hence, the difference between $\alpha$-invariant states and tracial states on $L(\Gamma)$ is similar to the different between inner amenability of $\Gamma$ and property Gamma of $L(\Gamma)$ (see e.g.\ \cite[p.394]{Vaes}).
It is shown in Proposition \ref{prop:st-inner}(a) that if $\Gamma$ is not strongly inner amenable, all $\alpha$-invariant states on $L(\Gamma)$ are tracial states.
On the other hand, if $\Gamma$ is an amenable ICC group, then $L(\Gamma)$ has only one tracial state, but there exist more than one $\alpha$-invariant states on $L(\Gamma)$ (by Theorem \ref{thm:ICC}(a)).
\end{rem}

\section{Invariant States and Property ($T$) (Theorem \ref{thm:property-T})}

Although Statement (G4) does not imply Statement (G2) for individual action, if Statement (G4) holds for every action, then so does Statement (G2).
This can be seen in the following result, which is a more elaborated version of Theorem \ref{thm:property-T}.

\begin{thm}\label{thm:equiv-prop-T}
Let $\Gamma$ be a countable discrete group.
The following statements are equivalent.

\begin{enumerate}[(T1)]
\item $\Gamma$ has property $(T)$.
\item If $\alpha$ is an action of $\Gamma$ on a von Neumann algebra $N$ and $\psi\in \CL(\KH)_+^*$ satisfying $\psi(\Cu_{\alpha_t}) = 1$ ($t\in \Gamma$), one has $\psi(p_{\Cu_{\alpha}}) \neq 0$.
\item $\KM_{N, \alpha} = \KM^{N, \alpha}$ for every action $\alpha$ of $\Gamma$ on a von Neumann algebra $N$.
\item For each action $\alpha$ of $\Gamma$ on a von Neumann algebra $N$ with $\dim N_*^\alpha  =1$, there is only one $\alpha$-invariant state on $N$.
\item For any Hilbert space $K$ and any action $\alpha$ of $\Gamma$ with $\CL(K)_*^\alpha  = (0)$, there is no $\alpha$-invariant state on $\CL(K)$.
\end{enumerate}
\end{thm}
\begin{prf}
Theorem \ref{thm:sp-gap-act-inv-mean} gives $(T1) \Rightarrow (T2) \Rightarrow (T3)$.
Moreover, it is clearly that $(T3) \Rightarrow (T4)$ and $(T3) \Rightarrow (T5)$. We need $(T4) \Rightarrow (T1)$ and $(T5) \Rightarrow (T1)$. 

\smnoind
$(T4) \Rightarrow (T1)$.
This follows from \cite[Theorem 2.5]{Sch}  (by considering the case when $N= L^\infty(X,\mu)$, where $(X,\mu)$ is a fixed non-atomic standard probability space, and $\alpha$ runs through all ergodic actions of $\Gamma$ on $(X,\mu)$).

\smnoind
$(T5) \Rightarrow (T1)$. 
Suppose that $\Gamma$ does not have property $(T)$.
By \cite[Theorem 1]{BV}, there is a unitary representation $\mu: \Gamma\to \CL(K)$ such that $\mu$ does not have a nonzero finite dimensional subrepresentation and $\pi:=\mu\otimes \overline{\mu}$ weakly contains $1_\Gamma$. 
It is easy to see, by using \cite[Proposition A.1.12]{BHV}, that $\pi$ does not have a nonzero finite dimensional subrepresentation as well. 

Let $N := \CL(K\otimes \overline{K})$ and $\alpha := {\rm Ad}\ \!\pi$. 
Then $\KH = K\otimes \overline{K}\otimes K\otimes \overline{K} = \CH\CS(K\otimes \overline{K})$ (where $\CH\CS$ denote the space of all Hilbert-Schmidt operators), $\Cu_{\alpha} = \pi \otimes \overline{\pi}$ and the $^*$-representation of $N$ on $\KH$ is given by compositions.
As $\pi$ does not have any nonzero finite dimensional subrepresentation, we see that $\KH^{\Cu_{\alpha}} = (0)$ (again by \cite[Proposition A.1.12]{BHV}).
On the other hand, if $\{\xi_i\}_{i\in I}$ is an almost $\pi$-invariant unit vector in $K\otimes \overline{K}$, the $\sigma(N^*,N)$-limit of a subnet of $\{\omega_{\xi_i}\}_{i\in I}$ will produce an $\alpha$-invariant state on $N$.
This contradicts Statement (T5).
\end{prf}

\begin{rem}
One may regard Theorem \ref{thm:equiv-prop-T} as a sort of non-commutative analogue of \cite[Theorem 2.5]{Sch}. 
There is also a non-commutative analogue of \cite[Theorem 2.4]{Sch}, which is basically a reformulation of \cite[Theorem 2.2]{B}. 
More precisely, the following are equivalence for a discrete group $\Gamma$:
\begin{enumerate}[({A}1)]
\item $\Gamma$ is amenable. 
\item $\KM_{N, \alpha} \neq \emptyset$ for every action $\alpha$ of $\Gamma$ on a von Neumann algebra $N$.
\item For any action $\alpha$ of $\Gamma$ on $\CL(K)$ (where $K$ is a Hilbert space) with $\CL(K)_*^\alpha = (0)$, there exists an $\alpha$-invariant state on $\CL(K)$. 
\end{enumerate}
In fact, if $\Gamma$ is amenable, then \cite[Theorem 5.1]{B} implies that $\Cu_\alpha\otimes \overline{\Cu_\alpha}$ weakly contains $1_\Gamma$. 
Thus, if we identify $N$ with the subalgebra $N\otimes 1$ of $\CL(\KH\otimes \KH)$, an almost $\Cu_\alpha\otimes \overline{\Cu_\alpha}$-invariant unit vector will produce an $\alpha$-invariant state on $N$. 
To show (A3) $\Rightarrow$ (A1), one may consider the action $\alpha := \ad \lambda$ on $\CL(\ell^2(\Gamma))$. 
It is easy to see that $\CL(\ell^2(\Gamma))^\alpha_* = (0)$ (e.g.\ by Lemma \ref{lem:H-N_*}(a)) and an $\alpha$-invariant state will restricts to a left invariant state on $\ell^\infty(\Gamma)$, which implies that $\Gamma$ is amenable. 
\end{rem}

Recall from \cite{LZ89} that $\Gamma$ has property $(T,FD)$ if $1_\Gamma$ is not in the closure of the subset of $\hat \Gamma$ consisting of non-trivial finite dimensional irreducible representations.
\begin{prop}\label{prop:no-inv-state}
Let $\Gamma$ be an infinite discrete group with property $(T,FD)$. 
Then (T1) is equivalent to the following statement. 
\begin{enumerate}[(T5')]
\item For any action $\alpha$ of $\Gamma$ on a von Neumann algebra $N$ with $N_*^\alpha = (0)$, there is no $\alpha$-invariant state on $N$.
\end{enumerate}
\end{prop}
\begin{prf}
By Theorem \ref{thm:sp-gap-act-inv-mean}, we have (T1)$\Rightarrow$(T5'). 
Now, suppose that $\Gamma$ does not have property $(T)$.
As $\Gamma$ has property $(T,FD)$, there is a net $\{(\pi^i,K^i)\}_{i\in I}$ in $\hat \Gamma\setminus\{1_\Gamma\}$ with each $K^i$ being infinite dimensional and there is a unit vector $\xi_i\in K^i$ ($i\in I$) such that
$$\|\pi_t^i\xi_i - \xi_i\|\ \to\ 0
\qquad (t\in \Gamma).$$
Let $N$ be the von Neumann algebra $\bigoplus_{i\in I} \CL(K^i)$ and set $\pi:=\bigoplus_{i\in I}\pi^i$ as well as $\alpha := {\rm Ad}\ \!\pi$.
Then $\KH = \bigoplus_{i\in I} \CH\CS(K^i)\cong \bigoplus_{i\in I} K^i\otimes \overline{K^i}$ and the representation of $N$ on $\KH$ is given by compositions.
In this case, one has
$$\Cu_{{\rm Ad}w}\big((\zeta^i \otimes \overline{\eta^i}\big)_{i\in I}\big)\ =\ \big(w_i\zeta^i \otimes \overline{w_i\eta^i}\big)_{i\in I}
\qquad \Big(w = (w_i)_{i\in I}\in {\bigoplus}_{i\in I}U(K^i); (\zeta^i \otimes \overline{\eta^i})_{i\in I}\in \KH \Big),$$
where $U(K^i)$ is the group of unitaries  on $K^i$.
Thus, $\Cu_{\alpha_t}(y) = \pi_t \circ y \circ \pi_{t^{-1}}$ ($y\in \KH; t\in \Gamma$), which implies that
$$\KH^{\Cu_{\alpha}}\ \subseteq\ {\bigoplus}_{i\in I}\{S_i\in \CH\CS(K^i): S_i\pi^i_t = \pi^i_t S_i,\ \! \forall t\in \Gamma\}$$
and hence $\KH^{\Cu_{\alpha}} = (0)$.

On the other hand, the $\sigma(N^*,N)$-limit of a subnet of  $\{\omega_{\xi_i}\}_{i\in I}$ (we consider $\CL(K^i)\subseteq N$ for all $i\in I$) will give an $\alpha$-invariant state on $N$.
\end{prf}

A similar statement as the above for the strong property $(T)$ of locally compact groups can be found in \cite{LN}. 

\begin{cor} 
Let $\Gamma$ be a minimally almost periodic group in the sense of \cite{vNW40} (i.e.\ there is no non-trivial finite dimensional irreducible representation of $\Gamma$).
If $\Gamma$ satisfies (T5'), then it has property (T) and hence is finitely generated.
\end{cor}

\section*{Acknowledgement}
We would like to thank the referee for helpful comments and suggestions that lead to a better presentation of the article, and for informing us that one can use the materials in \cite{BV} to obtain (T5) $\Rightarrow$ (T1) for general countable groups (we only had Proposition \ref{prop:no-inv-state} in the original version.)
We would also like to thank Prof. Roger Howe for helpful discussion leading to a simpler argument for Theorem \ref{thm:equiv-prop-T}. 

\bibliographystyle{plain}

\begin{thebibliography}{ZZZZ}

\bibitem
{B}
M.B. Bekka, Amenable unitary representations of locally compact groups, Invent.\ math.\ \textbf{100} (1990), 383-401.


\bibitem
{BC}
M.B. Bekka and Y. Cornulier, A spectral gap property for subgroups of finite covolume in Lie groups, Colloq.\ Math.\ \textbf{118} (2010), 175-182.


\bibitem
{BHV} M.B. Bekka, P. de la Harpe and A. Valette, \emph{Kazhdan's Property (T)}, Cambridge (2008).

\bibitem
{BV}
M.B. Bekka and A. Valette, Kazhdan's property (T) and amenable representations, Math.\ Z.\ \textbf{212} (1993), 293-299.


\bibitem
{Con}
A. Connes, Classification of injective factors. Cases $I\!I_1$, $I\!I_\infty$, $I\!I\!I_\lambda$, $\lambda \neq 1$, Ann.\ of Math.\
\textbf{104} (1976), 73-115.

\bibitem
{CW}
A. Connes and B. Weiss, Property T and almost invariant sequences,  Israel J.\ Math.\ \textbf{37} (1980), 209-210.

\bibitem
{Drin}
V. Drinfel'd, Finitely-additive measures on $S^2$ and $S^3$, invariant with respect to rotations, Funktsional.\ Anal.\ i Prilozhen \textbf{18} (1984), 77.

\bibitem
{Eff}
E.G. Effros, Property $\Gamma$ and inner amenability, Proc.\ Amer.\ Math.\ Soc.\ \textbf{47} (1975), 483-486.

\bibitem
{Green}
F. Greenleaf, \emph{Invariant means on topological groups and their applications}, Van Nostrand, Princeton, N. J. (1969).

\bibitem
{Haag-st-form}
U. Haagerup, The standard form of von Neumann algebras, Math.\ Scand.\ \textbf{37} (1975), 271-283.

\bibitem
{HaV}
P. de la Harpe and A. Valette, La propri\'{e}t\'{e} (T) de Kazhdan pour les groupes localement compacts, \emph{Ast\'{e}risque} \textbf{175} (1989).

\bibitem
{dR}
A. del Junco and J. Rosenblatt, Counterexamples in ergodic theory and number theory, Math.\ Ann.\ \textbf{245} (1979), 185-197.

\bibitem
{KRII}
R.V. Kadison and J.R. Ringrose, \emph{Fundamentals of the theory of operator algebras II: Advanced theory}, Grad.\ Stud.\ Math.\ \textbf{16}, Amer.\ Math. Soc.\ (1997).

\bibitem
{Kan-Mark}
E. Kaniuth and A. Markfort, The conjugation representation and inner amenability of discrete groups, J.\ reine angew.\ Math.\ \textbf{432} (1992), 23-37.

\bibitem
{LP}
A.T.M. Lau and A.L.T. Paterson, Inner amenable locally
compact groups, Trans.\ Amer.\ Math.\ Soc.\ \textbf{325} (1991), 155-169.

\bibitem{LN}
C.W. Leung and C.K. Ng, property $(T)$ of group homomorphisms, in preparation. 

\bibitem
{LH}
V. Losert and H. Rindle, Almost invariant sets, Bull.\ London Math.\ Soc.\ \textbf{13} (1981), 145-148.

\bibitem
{Lub}
A. Lubotzky, \emph{Discrete groups, Expanding graphs and Invariant measures}, Progress in Math. \textbf{125}, Birkh\"{a}user (1994).

\bibitem{LZ89}
A. Lubotzky and R. Zimmer, Variants Kazhdan's property for subgroups of semismiple groups, Israel J.\ Math.\ \textbf{66} (1989), 289-299. 

\bibitem
{Marg}
G.A. Margulis, Some remarks on invariant means, Mona.\ fur Math.\ \textbf{90} (1980), 233-235.

\bibitem
{Pasch}
W.L. Paschke, Inner amenability and conjugation operators, Proc.\ Amer.\ Math.\ Soc.\ \textbf{71} (1978), 117-118.

\bibitem
{Pat}
A.L.T. Paterson, \emph{Amenability}, Math.\ Sur.\ Mono.\ \textbf{29}, Amer.\ Math.\ Soc.\ (1988).

\bibitem
{Sch}
K. Schmidt, Amenability, Kazhdan's property T, strong ergodicity and invariant means for ergodic group-actions, Ergod.\ Th.\ \& Dynam.\ Sys.\ \textbf{1} (1981), 223-236.

\bibitem
{Sull}
D. Sullivan, For $n>3$ there is only one finitely additive rotationally invariant measure on the $n$-sphere on all Lebesgue measurable sets, Bull.\ Amer.\ Math.\ Soc.\ \textbf{4} (1981),  121-123.

\bibitem
{Vaes}
S. Vaes, An inner amenable group whose von Neumann algebra does not have property Gamma, Acta Math.\ \textbf{208} (2012), 389-394.

\bibitem
{vNW40}
J. von Neumann and E.P. Wigner, Minimally almost periodic groups, Annals of Math.\ \textbf{41} (1940), 746-750.
\end{thebibliography}

\end{document}